\documentclass[11pt]{article}

\usepackage{latexsym,amssymb,amsmath,amscd,amscd,amsthm,amsxtra,epsfig,setspace,authblk,psfrag}
		
\usepackage[usenames]{color}

\newtheorem{thm}{Theorem}

\newtheorem{lem}{Lemma}

\newtheorem{rmk}{Remark}

\newenvironment{pf}{{\noindent \it \bf Proof.}}{{\hfill$\Box$}\\}

\def\p{\partial}

\def\({\left(}
\def\){\right)}

\author{Zhou Zhang\footnote{‎zhangou@maths.usyd.edu.au}}
\affil{School of Mathematics and Statistics, 
University of Sydney, NSW 2006, Australia}

\author{Xiaoming Zheng\footnote{zheng1x@cmich.edu.}}
\affil{Department of Mathematics, Central Michigan University,
              Mount Pleasant, MI 48859 }   
\date{}
\title{The Representation of Line Dirac Delta Function Along a Space Curve}

\begin{document}
\maketitle
\begin{abstract} 
In this paper, we derive a new formula of the line Dirac delta function along a  curve in three-dimensional space: 
$\delta_\Sigma=\frac{\delta(\rho)}{2\pi \rho}$, where $\delta$ is a one-dimensional Dirac delta function and $\rho$ is the distance function to the curve $\Sigma$. Its extensions to level set formulation and submanifolds of co-dimensions 2 and 1 are also developed.  The main ideas can be applied for general dimension and codimension.
\end{abstract}


\vspace{0.2in}

\section{Introduction}
A line source expressed with a Dirac delta function is often used in mechanical, electromagnetic, and biological problems when a long and thin structure is involved. For instance, a line Dirac delta function is developed to address a line load acting on a plane \cite{SoedelPowder1979}. A distribution of electric charge on a curve can be represented by a line source \cite{Bladel1996}. 
A flexible structure such as an elastic fiber or filament is modeled as a line source  to represent the mechanical interactions between the structure and biofluid \cite{Zhu2002, Peskin2002}.
A very thin blood capillary in a bulk tumor can be treated as a line source when modeling a growth factor's reaction and diffusion \cite{Zhenginpress}. In these models, the line source is considered as a simple smooth curve $\Sigma$, embedded in the three-dimensional (3-D) Euclidean space $\mathbb{R}^3$ or its generalization. 

Denote the arc length parameter by $s$ and the corresponding spatial point on this curve by ${\bf x}(s)$. The line delta function $\delta_{\Sigma}$ associated with the curve $\Sigma$ is defined as a distribution such that  for any test function $f({\bf x}) \in C_{c}^{\infty} (\mathbb{R}^3)$ (i.e., infinitely differentiable and compactly supported), 
\begin{equation}\label{delta_def}
\int_{\mathbb{R}^3} \delta_{\Sigma}({\bf x})\cdot f({\bf x}) d{\bf x}=\int_\Sigma f\({\bf x}(s)\) ds.
\end{equation}
In Cartesian coordinates,   by formally switching the order of integrations (i.e., applying ``Fubini's Theorem''), it can be written as 
({\em e.g.}, see\cite{Peskin2002}) 
\begin{equation}\label{oldform}
\delta_{\Sigma}({\bf x}) =  \int_\Sigma \delta^{3D}\({\bf x} - {\bf x}(s)\) ds, \forall {\bf x} \in \mathbb{R}^3,
\end{equation}
where $\delta^{3D}({\bf x})=\delta(x) \delta(y) \delta(z)$ for ${\bf x}=(x,y,z)$, where $\delta(x)$ is the one-dimensional (1-D) Dirac delta function satisfying $\int_{-\infty}^{\infty}\delta(t) g(t)dt=g(0)$ for any $g\in C^\infty_c(\mathbb{R})$. 

The purpose of  this work is to provide a simpler representation of the line delta function $\delta_\Sigma$ for a fairly general class of spatial curves which is sufficient for most practical applications. More precisely, 
\begin{equation}
\label{eq:want}
\delta_{\Sigma}({\bf x}) =  \frac{\delta(\rho({\bf x}))}{2\pi \rho({\bf x})},
\end{equation} 
where $\rho({\bf x})$ is the distance function to the curve. This $\delta$ function is under the convention that $\int_0^\infty\delta(t) g(t)dt=g(0)$ for any $g\in C^\infty_c(\mathbb{R})$. By all means, the integration is over the whole natural domain of the variable for the $\delta$ function.

The main challenge is  to properly formulate the calculation of  $\int_{\mathbb{R}^3}\frac{\delta(\rho)}{2\pi \rho} f({\bf x}) d{\bf x}$ for any $f\in C_{c}^{\infty} (\mathbb{R}^3)$. For this, one has to choose the most reasonable and practical coordinate system. Indeed, the apparently singular term $\frac{1}{\rho}$ on the right hand side shows up naturally from the local cylindrical coordinates in a tubular neighbourhood of the curve. 
If the curve $\Sigma$ is simple and smooth, then a tubular neighbourhood can be readily chosen where the distance function works as the radius and will be cancelled out in coordinate transformation. This will be done in Section 2. In general, if the curve has non-smooth points, we have to further clarify the meaning of $\frac{\delta(\rho)}{2\pi \rho}$, which will be done in Section 3. 

\begin{figure}[htbp]
\begin{center}
\psfig{file=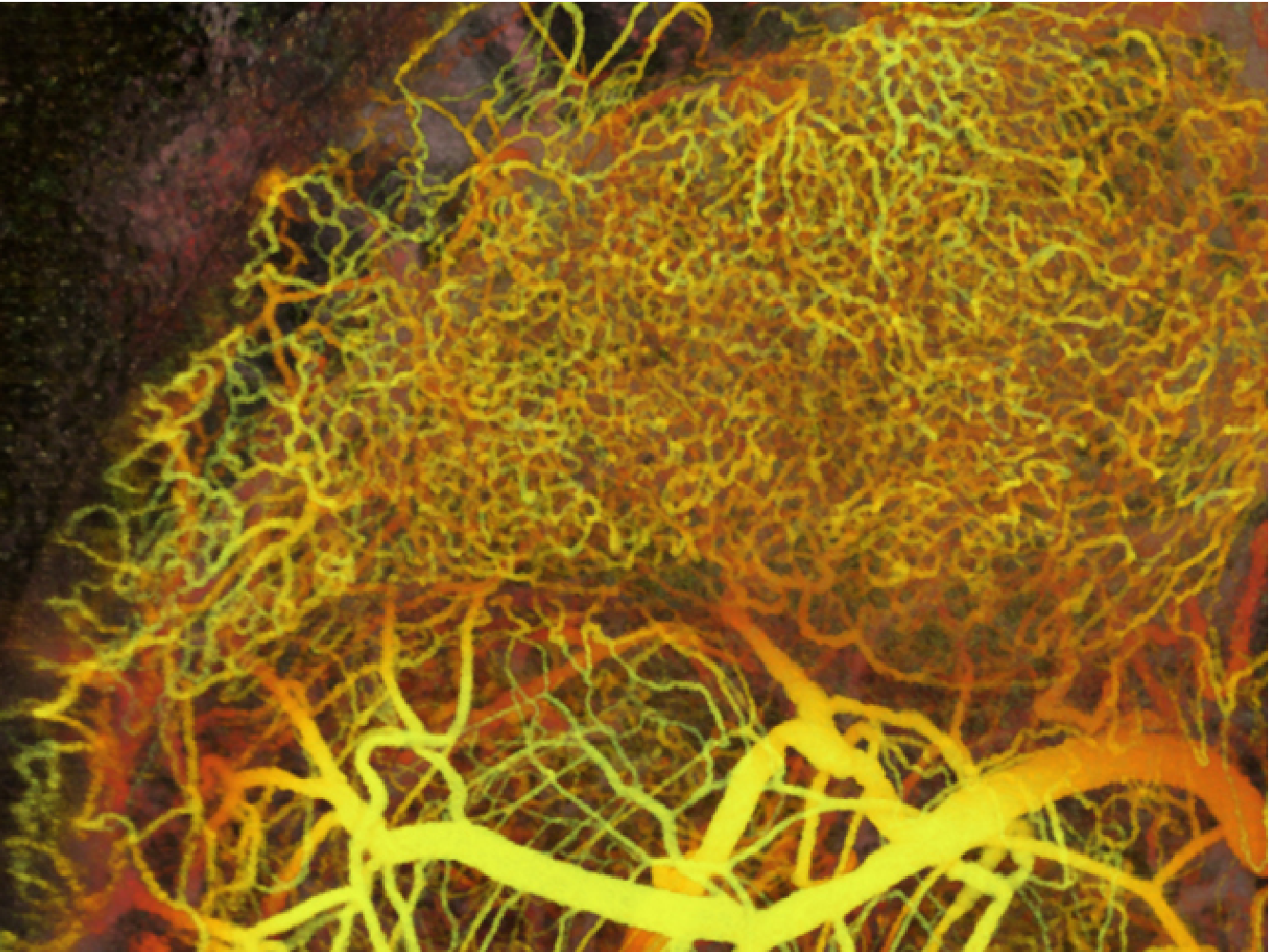, width=2.2in}
\caption{The highly irregular and tortuous blood vessel capillaries in a xenotransplanted U87 human glioblastoma multiforme tumor (upper part) in a mouse brain. The size of tissue shown in this figure is $2.6\,mm$ by $2\,mm$. This picture is taken from \cite{Vakoc2009} with permission.}
\label{tumor_vasculature}
\end{center}
\end{figure}

As an application, the formula \eqref{eq:want} works for  very complicated structure composed by piecewise smooth curves such as the newly formed blood capillary network in growing tumor, as shown in Fig.\,\ref{tumor_vasculature}. 
Another important application of this work is in the level set method.  The delta function of a curve or a surface represented by a level set function is crucial in the level set method \cite{Osherbook2003}, but the form for a spatial curve is not previously available. Eq.\, \eqref{eq:want} provides the first representation with the distance function, whose generalization to an arbitrary level set function will be given in Section 4.  The variation to plane curves in $\mathbb{R}^2$ is also discussed in Section 4. Finally, in Section 5, we further describe the extension to any dimension and co-dimension, and even manifolds. 

\section{Simple Smooth Curve\label{case1}}
The essential tool in our calculation of $\frac{\delta(\rho)}{2\pi \rho}$ is the knowledge of a neighbourhood around the curve. The standard tubular neighbourhood theorem is for compact smooth submanifolds of $\mathbb{R}^k$ (c.f. \cite{Bredon1993}). Therefore, when it is applied to a curve, the curve needs to be a submanifold of $\mathbb{R}^3$ and is diffeomorphic to a circle. However, because of the local nature of our problem, it can be readily applied to any simple smooth curve which is a submanifold of $\mathbb{R}^3$. More detail is as follows. 

Let $\Sigma$ be a simple smooth curve with the arc length parameter $s$ and any point on the curve denoted as ${\bf x}(s)$. Look at a finite piece of it with $s\in [s_1,s_2]$, for simplicity, still denoted as $\Sigma$. It is a closed set in $\mathbb{R}^3$ with the two endpoints ${\bf x}(s_1)$ and ${\bf x}(s_2)$.  Denote $\mathring{\Sigma}$ as the interior part by removing the two endpoints of the curve $\Sigma$.    
On each point ${\bf x}(s)\in\mathring{\Sigma}$, denote the tangent space as $T_{{\bf x}(s)}\Sigma$ which is exerted by $\frac{\partial}{\partial s}$. Define the normal plane as $N_{{\bf x}(s)}\Sigma=\{ {\bf z}\in \mathbb{R}^3 | {\bf z}\perp T_{{\bf x}(s_1)}\Sigma  \}$. On the two endpoints ${\bf x}(s_1)$ and ${\bf x}(s_2)$, define tangent and normal planes as the limits of those of points in $\mathring{\Sigma}$. Let $B({\bf x}(s_1),\epsilon)$ be the open ball centered at ${\bf x}(s_1)$ with radius $\epsilon$, and $B^{out}({\bf x}(s_1),\epsilon)$ be the open semi-ball of $B({\bf x}(s_1),\epsilon)$ cut by the normal plane $N_{{\bf x}(s_1)}\Sigma$ which is disjoint with $\mathring{\Sigma}$. Similarly we can define $B^{out}({\bf x}(s_2),\epsilon)$.

\begin{figure}[htbp]
\begin{center}
\psfrag{U}{$U$}
\psfrag{V}{$V$}
\psfrag{ds}{$\frac{\partial }{\partial s}$}
\psfrag{A}{$ $}
\psfrag{B}{$\rho({\bf x})$}
\psfrag{C}{${\bf x}(u,v,s)$}
\psfrag{D}{$\epsilon$}
\psfrag{E}{$\Sigma$}
\psfrag{F}{${\text{Tub}}^{\epsilon}\Sigma$}
\psfrag{G1}{$N_{{\bf x}(s_1)}\Sigma$}
\psfrag{G2}{$N_{{\bf x}(s_2)}\Sigma$}
\psfig{file=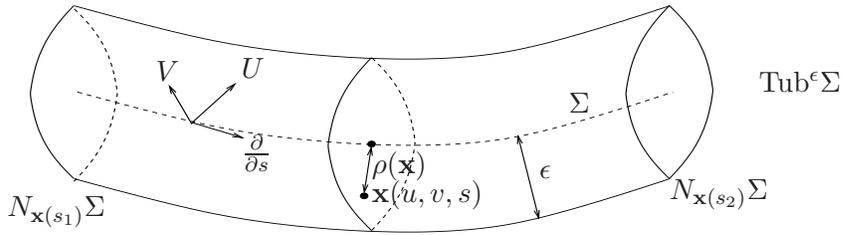, width=4in, height=1.2in}
\caption{A tubular neighbourhood ${\text{Tub}}^{\epsilon}\Sigma$ of a simple smooth curve $\Sigma$ in 3-D space bounded by two normal planes $N_{{\bf x}(s_1)}\Sigma$ and $N_{{\bf x}(s_2)}\Sigma$.
}
\label{fig2}
\end{center}
\end{figure}

Define the tubular neighbourhood 
$$\text{Tub}^{\epsilon}(\Sigma)=\{ {\bf x}\in\mathbb{R}^3|  dist(x,\Sigma)<\epsilon \} \backslash \{ B^{out}({\bf x}(s_1),\epsilon) \cup B^{out}({\bf x}(s_2),\epsilon) \},$$
which is the set of all points within distance $\epsilon$ to $\Sigma$ but excluding the two outside semi-balls at the two ends of the curve. Note  $N_{{\bf x}(s_1)}\Sigma$ and $N_{{\bf x}(s_2)}\Sigma$ are two bounding normal planes of this tube (see Fig.\,\ref{fig2}). Also define 
$$\Xi(\Sigma, \epsilon)=\{ ({\bf x},{\bf z}) \in \Sigma\times \mathbb{R}^3 | {\bf z}\in N_{{\bf x}(s)}\Sigma, ||{\bf z}||< \epsilon \}.$$
Let $\theta:\Xi(\Sigma, \epsilon)\to \mathbb{R}^3$ be given by $\theta({\bf x},{\bf z})={\bf x}+{\bf z}$.  The following is a standard result in differential geometry. 

\begin{lem}\label{lem1}
There exists $\epsilon>0$ such that the map $\theta$ is a diffeomorphism from $\Xi({\Sigma}, \epsilon)$ onto  $\text{Tub}^{\epsilon}\Sigma$.
\end{lem}

\begin{rmk}\label{rmk1}
A prominent feature of this tubular neighbourhood is that $\forall {\bf y}\in \text{Tub}^{\epsilon}\Sigma$, 
$dist({\bf y},\Sigma)=||z||$ where the vector ${\bf z}$ is from the unique decomposition ${\bf y}={\bf x} + {\bf z}$, $({\bf x},{\bf z})\in\Xi({\Sigma}, \epsilon)$. 
Also notice that for our curve of infinite length, one might not have a uniform $\epsilon$ for the whole curve, but this is not going to cause any trouble in this work because only pieces of finite length will be considered in the proof of Theorem \ref{thm1}.
\end{rmk}
Define the distance function to $\Sigma$ as $\rho({\bf x}) \triangleq\inf_{{\bf y}\in \Sigma}d({\bf x}, {\bf y})$, where $d({\bf x}, {\bf y})$ is the standard Euclidean distance. $\rho({\bf x})=0$ implies ${\bf x}\in\Sigma$ because the curve $\Sigma$  under consideration is always a closed subset of $\mathbb{R}^3$.
\begin{thm}\label{thm1}
Let $\Sigma$ be a simple smooth curve embedded in $\mathbb{R}^3$, then
\begin{equation}
\label{mainthm}
\delta_{\Sigma}({\bf x}) =  \frac{\delta(\rho)}{2\pi \rho}.
\end{equation}
\end{thm}
{{\noindent \it \bf Proof.}}
Around any point on $\Sigma$, we have a tubular neighbourhood $\text{Tub}^{\epsilon}\Sigma$ satisfying the condition in Lemma \ref{lem1}. 
First, a coordinate system will be constructed in $\text{Tub}^{\epsilon}\Sigma$ as follows. On the curve $\Sigma$, $\frac{\p}{\p s}$ is a unit vector field. Because the curve is smooth, there exist two other smooth unit vector fields $U$ and $V$ such that  $U$, $V$, and  $\frac{\p}{\p s}$ form an orthonormal system and satisfy the right hand rule. On each normal plane $N_{{\bf x}(s)}\Sigma$, $U$ and $V$ provide a $\mathbb{R}^2$ coordinate system, $(u, v)$. Denote any point on this plane as ${\bf x}(u, v, s)$ (see Fig.\,\ref{fig2}).
According to Remark \ref{rmk1}, $\rho\bigl({\bf x}(u, v, s)\bigr)=\sqrt{u^2+v^2}$. 

Because $\{U, V, \frac{\partial}{\partial s}\}$ is a smooth frame along the curve, the map from $(u, v, s)$ to $(x, y, z)$ is smooth.
The Jacobian between $(x, y, z)$ and $(u, v, s)$ is  the determinant of the transition matrix between the $\mathbb{R}^3$ standard basis $\{\frac{\partial}{\partial x}, \frac{\partial}{\partial y}, \frac{\partial}{\partial z}\}$ and $\{\frac{\partial}{\partial u}, \frac{\partial}{\partial v}, \frac{\partial}
{\partial s}\}$ which is just $\{U, V, \frac{\p}{\p s}\}$ when restricted to the curve.
Although $ \frac{\partial}{\partial s}$ is not explicit away from the curve, the transition matrix has to be in $SO(3)$  on the curve because it is between two orthonormal bases with the same (right-handed) orientation. So the determinant is $1$ on the curve. As the Jacobian is smooth in the tubular neighbourhood $\text{Tub}^\epsilon\Sigma$,  it's $1+O(\sqrt{u^2+v^2})$ in light of the Taylor series, where we always only care about the local behavior for some finite piece of the curve.   
Thus the Jacobian must be nonzero in a sufficiently small neighbourhood, say, $\text{Tub}^{\epsilon_1}\Sigma,0<\epsilon_1<\epsilon$. For simplicity, it is still denoted as $\text{Tub}^\epsilon\Sigma$. Therefore, $(u, v, s)$ is a coordinate system in $\text{Tub}^\epsilon\Sigma$.  

When converted to the polar version $(\rho, \theta, s)$, since $\det \left(  \frac{\partial (u, v, s)}{\partial (\rho,\theta,s)}\right)=\rho$, it is easy to see 
$$\det \left(  \frac{\partial (x,y,z)}{\partial (\rho,\theta,s)}\right)=\rho+O(\rho^2) ~\text{in}~\text{Tub}^\epsilon\Sigma.$$
 
For any test function $f({\bf x}) \in C_{c}^{\infty} (\mathbb{R}^3)$,  $\int_{\mathbb{R}^3}\frac{\delta(\rho)}{2\pi \rho} f({\bf x}) d{\bf x}$ only depends on the situation near $\Sigma$ as $\delta(\rho)$ is supported on $\Sigma$. Since $f$ is compactly supported, we can assume that it is supported in a ball $B$. Because $B\cap \Sigma$ is compact, there exist finitely many local tubular neighbourhoods centered at involved pieces of $\Sigma$  covering $B\cap \Sigma$,  denoted as $\{T_i\}_{i\in I}$. We can also make sure that $\{\widehat T_i\}$ with $\widehat T_i\Subset T_i$ also covers $B\cap \Sigma$. Now a stardard partition of unity construction provides us with $h_i\in C^\infty_c(T_i)$ such that $\sum_{i\in I}h_i=1$ when restricted to $\cup_{i\in I}\widehat T_i$. Finally, we can carry out the following computation. 
\begin{eqnarray}
\int_{\mathbb{R}^3}\frac{\delta(\rho)}{2\pi \rho}  f({\bf x}) d{\bf x}
\nonumber 
&=& \int_{\mathbb{R}^3}\sum_{i\in I}\frac{\delta(\rho)}{2\pi \rho}  f({\bf x})h_i({\bf x}) d{\bf x} 
=  \sum_{i\in I} \int_{\widehat T_i} \frac{\delta(\rho)}{2\pi \rho}  f({\bf x})h_i({\bf x}) d{\bf x}   
\nonumber\\
&=& \sum_{i\in I}\int_s\int_\theta\int_\rho \frac{\delta(\rho)}{2\pi \rho} f(\rho,\theta,s) h_i(\rho,\theta,s) \(\rho+O(\rho^2)\) d\rho d\theta ds \nonumber \\
&=&  \sum_{i\in I}\int_s\int_\theta\int_\rho  \frac{\delta(\rho)}{2\pi} f(\rho,\theta,s) h_i(\rho,\theta,s)\(1+O(\rho)\) d\rho d\theta ds \nonumber \\
&=& \sum_{i\in I}\int_s\int_0^{2\pi} \frac{1}{2\pi} f(0,\theta,s) h_i(0,\theta,s) d\theta ds\nonumber \\
&=& \int_s f({\bf x}(s)) \sum_{i\in I} h_i({\bf x}(s)) ds 
= \int_s f({\bf x}(s)) ds =\int_\Sigma f({\bf x}(s)) ds. \nonumber
\label{feifei5}
\end{eqnarray}
In light of the definition in Eq.\,\eqref{delta_def}, the theorem is proven. 

{{\hfill$\Box$}\\} 

\begin{rmk}
\label{rmk2}
The choice of coordinates does not have to be fixed and is actually flexible. However, we do need to start with a Euclidean orthonormal basis along $\Sigma$ with one vector in the direction of the curve. In short, the delta function would be the function on the normal plane of the curve with the Euclidean structure induced from $\mathbb{R}^3$,  which is of course the most natural point of view. 

\end{rmk}

\section{General Curve  \label{case2}}
With a proper understanding of the delta function, the previous argument can be adjusted to a general class of curves which are `` topological graphs with smooth edges". The graph can have infinite vertices with finite edges between any two of them and the edges need not be straight (i.e., being topological). We only need to require that there is no local ``accumulation'' happening, which is in general what people think of graphs.
  
\begin{figure}[htbp]
\begin{center}
\psfig{file=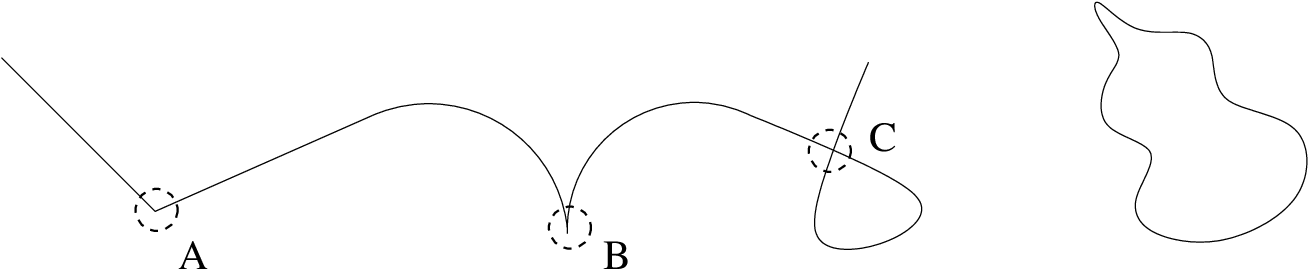, width=5in, height=1in}
\caption{This curve has two connected components: the left component is non-closed, has non-smooth points at A and B, and crosses itself at C; the right component is closed and smooth.}
\label{allcurves}
\end{center}
\end{figure}

Let $\Sigma$ be such a curve. The whole curve $\Sigma$ may be closed or open, may have non-smooth points and self-intersections, and may have more than one connected components (See Fig.\,\ref{allcurves}). 
Denote the collection of all non-smooth points as $P_\Sigma=\{{\bf p}_j\}_{j\in \Gamma}$, where $\Gamma$ is a finite or countable index set. To handle these non-smooth points, we remove from $\Sigma$ a small neighbourhood in $\mathbb{R}^3$ around each such point such that the boundary of the neighbourhood is perpendicular to the curve at the point of intersection.  Denote the union of these neighbourhoods by $U_\epsilon$ and shrink $U_\epsilon$ to $P_\Sigma$ as $\epsilon\to 0$. Let $\chi_\epsilon$ be the characteristic function of $\mathbb{R}^3\setminus U_\epsilon$, i.e., $\chi_\epsilon=0$ over $U_\epsilon$ and $1$ otherwise. 

The following definition is what we want: 
\begin{equation}
\label{def:delta}
\frac{\delta(\rho)}{2\pi \rho} = \lim\limits_{\epsilon\to 0}\frac{\chi_\epsilon\cdot\delta(\rho)}{2\pi \rho}.
\end{equation}
Notice that $\frac{\chi_\epsilon\cdot\delta(\rho)}{2\pi \rho}$ can be understood in the same way as before because a neighbourhood of the non-smooth points has been removed. The limit is in the weak sense, i.e., in the sense of current (or integration), with the  existence justified by the following computation. For any $f\in C^\infty_c(\mathbb{R}^3)$,  
\begin{eqnarray}
\lim\limits_{\epsilon\to 0} \int_{\mathbb{R}^3}\frac{\chi_\epsilon\cdot\delta(\rho)}{2\pi \rho}  f({\bf x}) d{\bf x} \nonumber 
= \lim\limits_{\epsilon\to 0} 
\int_{\Sigma} \chi_\epsilon\cdot f({\bf x}(s)) ds 
=  \int_{\Sigma} f({\bf x}(s)) ds\nonumber
\end{eqnarray}
The first equality makes use of the same kind of computation as that in the previous section, i.e., using a partition of unity construction to reduce to local tubular neighbourhoods with cylindrical coordinates. This is done for each $\epsilon$ and the local tubular neighbourhoods only need to cover $\Sigma\setminus U_\epsilon$.  

This also shows the independence of the choice (on $\chi_\epsilon$) in the construction. Simply speaking, there is no contribution from $P_\Sigma$ in integration. In summary, we have proven the following main theorem.
\begin{thm}\label{thm2}
Let $\Sigma$ be a topological graph in $\mathbb{R}^3$ with smooth edges and $\rho$ be the distance function to $\Sigma$. Then the definition in Eq.\,\eqref{def:delta} is independent of choices and 
$\delta_{\Sigma}({\bf x}) =  \frac{\delta(\rho)}{2\pi \rho}$. 
\end{thm}

\section{Level Set Method \label{level-set}}

Still consider in $\mathbb{R}^3$. Let $\phi$ be a non-negative level set function of the curve $\Sigma$, i.e., $\phi\geqslant 0$ and $\Sigma=\{{\bf x}~|~ \phi({\bf x})=0 \}$. $\phi$ is a function, $\phi=\phi(\rho, \theta, s)$, locally near $\Sigma$ in the coordinate system used before, and assume $\frac{\p\phi}{\p\rho}\ne 0$ around $\Sigma$. Since $\phi$ is non-negative, we have actually assumed $\frac{\p\phi}{\p\rho}>0$. 
This implies $(\phi, \theta, s)$ can be viewed as a generalized cylindrical coordinate system locally near $\Sigma$, and the Jacobian between it and $(\rho, \theta, s)$ is $\frac{\p(\phi, \theta, s)}{\p(\rho, \theta,s)}=\frac{\p\phi}{\p\rho}$. Then, we have the following result. 
\begin{thm}
Let $\Sigma$ be a topological graph in $\mathbb{R}^3$ with smooth edges, $\rho$ be the distance function to $\Sigma$, and $\phi({\bf x})$ be a non-negative level set function of $\Sigma$ described above, then 
$\delta_{\Sigma}({\bf x}) =  \frac{\delta(\phi)\frac{\p\phi}{\p\rho}}{2\pi \rho}$.
\end{thm}

For the proof, one only needs to notice $d\rho d\theta ds=\frac{1}{\p\phi/\p\rho}d\phi d\theta ds$. Non-smooth points can be treated in the same way as discussed in Section 3.

\begin{rmk}

Let's point out that the above $\frac{\p \phi}{\p \rho}$ is not well-defined as a function over $\Sigma$ in general because of the dependence on $\theta$, and so the conclusion of Theorem 3 should really be understood in the sense of current (or integration). 

\end{rmk} 

Using the formulation as for Theorem 3, we have the similar discussion for plane curves in $\mathbb{R}^2$. In $\mathbb{R}^2$, let $\phi$ be a non-negative level set function of the curve $\Sigma$, i.e., $\phi\geqslant 0$ and $\Sigma=\{{\bf x}~|~ \phi({\bf x})=0 \}$. Now the corresponding $\theta$ takes value in $S^0=\{1, -1\}$ instead of $S^1$, and $\frac{\p\phi}{\p\rho}$ has two values for each point on $\Sigma$ depending on the choice of $\theta$ value. If we assume them to be positive, then we end up with $\delta_{\Sigma}({\bf x}) =\frac{1}{2}\delta(\rho)= \frac{1}{2}\delta(\phi)\frac{\p\phi}{\p\rho}$.    

By all means, this is not very natural for the plane curve setting and we have the following alternative and better formulation for plane curve. 
In $\mathbb{R}^2$, let $\phi$ be a level set function of the curve $\Sigma$, i.e., $\Sigma=\{{\bf x}~|~ \phi({\bf x})=0 \}$. We assume $|\nabla \phi|\ne 0$ along $\Sigma$. Now we define the ``signed'' distance function to $\Sigma$, denoted by $\tilde\rho$, by giving the same sign as the $\phi$ value for the point under consideration. This $\tilde\rho$ behaves nicely near $\Sigma$ which is all we need. 

As $|\nabla\phi|\neq 0$ over $\Sigma$, we have $\frac{\p \phi}{\p \tilde\rho}>0$ along (and so near) $\Sigma$. We then have the following theorem. Of course, the corresponding delta function $\delta(\tilde\rho)$ will be integrated over $\mathbb{R}$.   

\begin{thm}
Let $\Sigma$ be a topological graph in $\mathbb{R}^2$ with smooth edges, $\tilde\rho$ be the signed distance function to $\Sigma$ and $\phi({\bf x})$ be a level set function of $\Sigma$ as described above. Then we have 
$$\delta_{\Sigma}({\bf x}) = \delta(\tilde\rho)= \delta(\phi) |\nabla\phi|.$$
\end{thm}
\begin{pf}
In light of the discussion in Section 3, we only need to prove this for a simple smooth curve $\Sigma$ embedded in $\mathbb{R}^2$. We use the 2-D coordinates $(\tilde\rho, s)$ in the tubular neighourhood $U$ around $\Sigma$. Over $\Sigma$, $\nabla\phi = \hat{\rho} \frac{\p\phi}{\p\tilde\rho} + \hat{s} \frac{\p\phi}{\p s}= \hat{\rho} \frac{\p\phi}{\p\tilde\rho}$ as $\frac{\p\phi}{\p s}=0$ over $\Sigma$, where $\hat{\rho}$ and $\hat{s}$ are orthonormal vectors in the direction of $\frac{\p}{\p\tilde\rho}$ and $\frac{\p}{\p s}$, respectively. Then over $\Sigma$, $|\nabla\phi|=\frac{\p\phi}{\p\tilde\rho}>0$.
Take any test function $f\in C^{\infty}_c(\mathbb{R}^2)$,
\begin{eqnarray*}
\int_{\mathbb{R}^2} f({\bf x}) \delta(\tilde\rho)d{\bf x}
&=&\int_U f({\bf x})  \delta(\tilde\rho)J d\tilde\rho ds \\
&=& \int_{\Sigma} f({\bf x}(s)) ds,
\end{eqnarray*}

\begin{eqnarray*}
\int_{\mathbb{R}^2} f({\bf x}) \delta(\phi) |\nabla\phi| d{\bf x}
&=&\int_U f({\bf x}) \delta(\phi) |\nabla\phi| J d\tilde\rho ds \\
&=& \int_U f({\bf x}) \delta(\phi)  |\nabla\phi| J \frac{1}{{\p\phi}/{\p\tilde\rho}} d\phi ds\\
 &=& \int_{\Sigma} f({\bf x}(s)) ds,
\end{eqnarray*}
where  $J$ is the Jacobian between the Euclidean coordinate and $(\tilde\rho, s)$, which is equal to $1$ along $\Sigma$ because $\frac{\p}{\p\tilde\rho}$ and $\frac{\p}{\p s}$ form an orthonormal basis there.

\end{pf}

\begin{rmk}

In  the 3-D case, a succinct formula like $\frac{\delta(\phi)}{2\pi \rho} |\nabla\phi|$ does not hold in general. Indeed, in 3-D cylindrical coordinates, $|\nabla\phi|\thicksim \sqrt{\phi_{\rho}^2 + \frac{1}{\rho^2} \phi_{\theta}^2 + \phi_s^2}$ near $\Sigma$. Therefore, when approaching $\Sigma$, $\frac{\p \phi}{\p \rho}$ can't be replaced by $|\nabla\phi|$ unless $\phi_{\theta}/\rho\to 0$, which puts serious restriction on this level set function $\phi$. 

\end{rmk}

\section{Further Extensions}

The study here can naturally be extended to more general case of a submanifold, $\Sigma$, of co-dimension 2 or 1 in $\mathbb{R}^n$ with $n>3$ (or any smooth Riemannian manifold). That is, $\delta_\Sigma=\frac{\delta(\rho)}{2\pi\rho}$ for a submanifold $\Sigma$ with co-dimension 2, while $\delta_\Sigma=\delta(\tilde\rho)=\delta(\phi) |\nabla\phi|$ if $\Sigma$ is of co-dimension 1 in the level set method setting. For even higher co-dimension case, one can proceed in exactly the same way as for the co-dimension 2 case, with $\theta$ taking value in some higher dimensional sphere.   

Indeed, the form $\delta_\Sigma=\delta(\phi) |\nabla\phi|$ for smooth closed curves in $\mathbb{R}^2$ and surfaces in $\mathbb{R}^3$ was first derived in \cite{Osher1996} and now is widely used in the level set method. The results presented in this work can be regarded as the extension of \cite{Osher1996}.

\section*{Acknowledgement}
Zhou Zhang thanks the University of Sydney and ARC DP110102654. Xiaoming Zheng thanks Central Michigan University ORSP Early Career Investigator grant \#C61373.

\end{document}